\documentclass[11pt]{article}

\usepackage{amsmath}
\usepackage{amssymb}
\usepackage{amsfonts}
\usepackage{amsthm}

\usepackage[utf8x]{inputenc} 

\usepackage{graphicx}
\usepackage[autostyle]{csquotes}
\usepackage[linktocpage]{hyperref}
\usepackage{fullpage}
\usepackage{slashed}
\usepackage{mathtools}
\usepackage{indentfirst}
\usepackage{bm}
\usepackage{mathrsfs} 
\usepackage{authblk}
\usepackage{tikz}
\usepackage{tikz-cd}

\usepackage{subfig}

\allowdisplaybreaks

\newtheorem{theorem}{Theorem}[section]
\newtheorem{prop}[theorem]{Proposition}
\newtheorem{cor}[theorem]{Corollary}
\newtheorem{lemma}[theorem]{Lemma}
\newtheorem{conjecture}[theorem]{Conjecture}

\newcommand{\intg}{\mathbb{Z}}
\newcommand{\rational}{\mathbb{Q}}

\newcommand{\rarw}{\rightarrow}

\newcommand{\lb}{\left(}
\newcommand{\rb}{\right)}
\newcommand{\lsb}{\left[}
\newcommand{\rsb}{\right]}
\newcommand{\lac}{\left\{}
\newcommand{\rac}{\right\}}

\newcommand{\RN}[1]{\textup{\uppercase\expandafter{\romannumeral#1}}}

\title{A Cable Knot and BPS-Series $\RN{2}$}
\author{John Chae}
\affil{University of California Davis, 1 Shields Ave, Davis, CA, 95616, USA \\ {yjchae@ucdavis.edu}}
\date{}

\begin{document}

\maketitle

\begin{abstract}
This is a companion paper to earlier work of the author, which generalizes to an infinite family of $(2,2w+1)$-cabling of the figure eight knot ($|w|>3$) and proposes general formulas for the two-variable series invariant of the family of the cable knots. The formulas provide an insight into the cabling operation. We verify the conjecture through explicit examples using the recursion method, which also provide a strong evidence for the $q$-holonomic property of the series invariant. This result paves a road for computation of the WRT invariant of a 3-manifold obtained from a Dehn surgery on the cable knots via a certain $q$-series. We also analyze and conjecture formulas for $(3,3w+1)$-cabling ($|w|>3$).

\end{abstract}


\section{Introduction}


Categorification of link invariants has been a fruitful program that has deepened and broadened our understanding of link invariants. Beginning with the advent of the
Khovanov homology that categorifies the (colored) Jones polynomials of links~\cite{Kh}(\cite{Kh2}), there has been a series of constructions of other link invariants such as categorification of $sl(N)$ polynomial by $sl(N)$ Khovanov-Rozannsky homology~\cite{KR} and the HOMFLY homology~\cite{KR2} that categorifies the HOMFLY polynomial. Some homological theories have more than one approaches for their computations. For example, in addition to the original matrix factorization formulation of \cite{KR2}, there are braid varieties~\cite{Trinh} and (local) Hilbert scheme of points~\cite{ORS} formulations\footnote{These other approaches only apply to particular classes of links.}. Another notable example is alternative categorification of the colored Jones polynomials as shown in \cite{BW, CK}. These examples illustrates multi-faceted character of the link homology theories tied to a variety of perspectives of links in 3-sphere. Inspired by these significant advancement for links, there has been an initiation of a categorification of a (numerical) invariant of closed oriented 3-manifolds called the Witten-Reshitikhin-Turaev (WRT) invariant~\cite{W,RT1} in \cite{Kh3, Q}.  Some years ago, a prediction for a categorification of the WRT invariant proposed from mathematical physics~\cite{GPPV}. It predicts a power series invariant of a closed oriented 3-manifold $Y$ often denoted by $\hat{Z} (q)$, which is conjecturally a graded Euler characteristic of a certain triply graded homology group denoted by $\mathcal{H}^{i,j}_{BPS}(Y;b)$ as follows
$$
\hat{Z}_{b}[Y;q]= \sum_{i,j} (-1)^i\, q^j\, dim\, \mathcal{H}^{i,j}_{BPS}(Y;b)\quad \in\, \frac{1}{2^c}\, q^{\Delta (Y)}\, \intg [[q]],\qquad \Delta \in \rational \qquad c \in \intg_{\geq 0} \qquad |q| < 1.
$$
This series invariant was extended to a complement $M^{3}_{K}$ of a knot $K \subset \intg HS^3$~\footnote{A generalization to links is given in \cite{P2}} in\cite{GM}, 
$$
F_{K}(x,q) : = \hat{Z} (M^{3}_{K}; x^{1/2},q)\qquad |q| < 1,
$$
which is a two variable series in which $x$ is associated to the relative $Spin^c (M^{3}_{K},T^2)$-structures. This knot invariant $F_K$ takes the form~\footnote{Implicitly, there is a choice of group; originally, the group used is ${\rm SU}(2)$.}
\begin{gather}
F_K(x,q)= \frac{1}{2} \sum_{\substack{m \geq 1 \\ m \ \text{odd}}}^{\infty} f_{m}(q) \big(x^{m/2}-x^{-m/2}\big) \in \frac{1}{2^{c}} q^{\Delta} \mathbb{Z}\big[x^{\pm 1/2}\big]\big[\big[q^{\pm 1}\big]\big],
\end{gather}
where $f_{m}(q)$ are Laurent series with integer coefficients\footnote{They can be polynomials for monic Alexander polynomial of $K$ (See Section 3.2)}, $c \in \mathbb{Z}_{\geq 0}$ and $\Delta \in \mathbb{Q}$. A curious feature of $F_K$ invariant is that it seems to have similar properties as that of the colored Jones polynomial. For example, it is conjectured to satisfy the Melvin–Morton–Rozansky theorem~\cite{MM,R,R2} (proven in~\cite{BG}): 
\begin{conjecture}[{\cite[Conjecture 1.5]{GM}}]\label{conjecture3} For a knot $K \subset$ $\intg HS^3$, the asymptotic expansion of the knot invariant $F_{K}\big(x,q={\rm e}^{\hbar}\big)$ about $\hbar =0$ coincides with the Melvin--Morton--Rozansky (MMR) expansion of the colored Jones polynomial in the large color limit:
\begin{gather}
\frac{F_{K}\big(x,q={\rm e}^{\hbar}\big)}{x^{1/2}-x^{-1/2}} = \sum_{r=0}^{\infty} \frac{P_{r}(x)}{\Delta_K(x)^{2r+1}}\hbar^r,
\end{gather}
where $x=e^{n\hbar}$ is fixed, n is the color of $K$, $P_{r}(x) \in \mathbb{Q} \big[x^{\pm 1}\big]$, $P_{0}(x)=1$ and $\Delta_K (x)$ is the (symmetrized) Alexander polynomial of~$K$.
\end{conjecture}
\noindent This property was proven for a wide class link in \cite{P}. Additionally, it possesses the q-holonomic property similar to that of the color Jones polynomials~\cite{GL}:
\begin{conjecture}[{\cite[Conjecture 1.6]{GM}}] For any knot $K \subset$ $\intg HS^3$, the normalized series $f_{K}(x,q)$ satisfies a linear recursion relation generated by the quantum A-polynomial $\hat{A}_K(q,\hat{x},\hat{y})$ of $K$ :
\begin{gather}
\hat{A}_{K}(q, \hat{x},\hat{y}) f_{K}(x,q) = 0,
\end{gather}
where $f_{K}:=F_{K}(x,q)/\big(x^{1/2}-x^{-1/2}\big)$. 
\end{conjecture}
\noindent The actions of $\hat{x}$ and $\hat{y}$ are
$$
\hat{x} f_{K}(x,q)= x f_{K}(x,q) \qquad \hat{y}f_{K}(x,q)= f_{K}(xq,q).
$$
\newline

In \cite{C}, $F_K (x,q)$ was computed for a cable knot of the figure eight knot $K=C_{9,2}(4_1)$ in $\intg HS^3$. Most importantly, a relation between $F_K (x,q)$ and $F_{4_1} (x,q)$ was found. In this paper, we generalize the relation to $(2w+1, 2)$ for $w>3$ and $(3w+1, 3)$ for $w>3$. To the best of the author's knowledge, the recursion method is the only available method to date for a computation of $F_K (x,q) $ for the cablings of the figure eight knot. The R-matrix formulation of $F_K (x,q) $ in \cite{P2} applies to certain classes for braids and hence links. Although they are broad, the cables of the figure eight considered in this paper may or may not fall into those classes. In any case, we propose a general formula for $F_K (x,q) $ for the above mentioned cablings of the figure eight knot.
\newline


\section{$(2w+1,2)$-cables\, $(w>3)$}

%

The quantum (or noncommutative) A-polynomial of a class of cable knot $K_r = C_{(r,2)}(\bm{4_1})$ in $S^3$ having minimal L-degree is given by~\cite{Ru,Tran}
\begin{equation}
\hat{A}_{K} (t,M,L)=(L-1)B(t,M)^{-1}Q(t,M,L) \lb M^r L + t^{-2r}M^{-r} \rb \in \tilde{\mathcal{A}}_{K},\quad gcd(r,2)=1, \quad |r| > 8
\end{equation}
The definitions of $Q(t,M,L)$ and $B(t,M)$ are written in Appendix A. Applying (4) to $f_K(x,q)$ together with $x=q^n$ yields
\begin{equation}
\alpha (x,q)F_{K_r}(x,q) + \beta (x,q)F_{K_r}(xq,q)+ \gamma (x,q)F_{K_r}(xq^2,q) + \delta (x,q) F_{K_r}(xq^3,q) + F_{K_r}(xq^4,q)=0,
\end{equation} 
where $\alpha ,\beta ,\gamma , \delta  $ are rational functions. We propose the following formula for initial conditions needed for (5) in terms of the coefficient functions $\lac h_s (q) \rac$ of  $F_{\bm{4_1}}(x,q)$. They provide insights into the cabling operation for $F_K (x,q)$ and exhibit simple patterns, which enable us to obtain the complete set of the initial conditions $\lac f_m (q) \rac $ for any value of $r$. We verify the initial conditions through explicit examples (see Section 3) and find that the patterns persist to $f_s (q)$ that are beyond the initial conditions. This provides a strong evidence that we can obtain $F_{C_{(r,2)}(\bm{4_1})} (x,q)$ to any desired order in $x^{\pm 1/2}$ efficiently without resorting to the laborious recursion method. 
\newline

\indent We state the nonzero initial conditions needed for (5) in the case of $r=2w+1,\, w>3$, which covers the most of the positive subspace of the cabling parameter space (see below for the negative case).

\begin{conjecture} 
For a $(2w+1,2)$-cabling of $\bm{4_1}\, (w > 3)$ in $\intg HS^3$, nonzero initial data $\lac f_m (q) \in \intg [q^{\pm 1}] \rac = \lac f^{+}_i (q)  |\, \text{all coefficients are positive} \rac \bigsqcup \lac f^{-}_j (q) |\, \text{all coefficients are negative}  \rac$ for (5) are determined by the coefficient functions $\lac h_s (q) \rac$ of  $F_{\bm{4_1}}(x,q)$ as follows.
\begin{align*}	
f^{+}_{2w+3} (q) & = 2\,  h_1 (q)\, q^n \qquad ( h_1 (q)=1)\\
f^{+}_{2w+7} (q) & = 2\, h_3 (q)\, q^{n+1}  \qquad ( h_3 (q)=2)\\
f^{+}_{2w+11} (q) & = 2\, h_5 (q)\, q^{n+2}\\
&\vdots\\
f^{+}_{10w+7} (q) & = 2\, \lb  h_{ 4w+3} (q)\, q^{n+2w+1 } + h_1 (q) \, q^{n +2w+1 +  \Delta}  \rb\\
f^{+}_{10w+11} (q) & = 2\, \lb h_{ 4w+5 } (q)\, q^{n+ 2w+2 } + h_3 (q) \, q^{n + 2w+6  + \Delta }  \rb\\ 
f^{+}_{10w+15} (q) & = 2\, \lb h_{ 4w+7 } (q)\, q^{n+  2w+3} + h_5 (q) \, q^{n + 2w+11  + \Delta }  \rb\\
&\vdots\\
f^{+}_{18w+11} (q) & = 2\, \lb h_{ 8w+5 } (q)\, q^{n+ 4w+2} + h_{4w+3 } (q) \, q^{n + 12w+6 +  \Delta }  + h_1 (q) \, q^{n + 12w+6 + 2\Delta }  \rb\\
f^{+}_{18w+15} (q) & = 2\, \lb h_{8w+7} (q)\, q^{n+ 4w+3 } + h_{4w+5} (q) \, q^{n +12w+11 +  \Delta }  + h_3 (q)\, q^{n + 12w+15 + 2\Delta }  \rb\\
&\vdots\\
f^{-}_{6w+5} (q) & = -2\,  h_1 (q)\, q^{3n} \\
f^{-}_{6w+9} (q) & = -2\, h_3 (q)\, q^{3n+3}  \\
f^{-}_{6w+13} (q) & = -2\, h_5 (q)\, q^{3n+6}\\
&\vdots\\
f^{-}_{14w+9} (q) & = -2 \lb h_{4w+3}(q) \, q^{3n+ 6w+3} + h_1 (q) \, q^{3n+ 6w+3 + \Delta} \rb\\
f^{-}_{14w+13} (q) & = -2 \lb h_{4w+5}(q) \, q^{3n+6w+6 } + h_3 (q) \, q^{3n+ 6w+10 + \Delta} \rb\\
&\vdots\\
\end{align*}
where $n \in \intg$ and $\Delta =4w+4$.
\end{conjecture}
\noindent We state observations about the above formulas.\\
\noindent \textbf{Remark 1} Since the order of the $\lac f_m (q) \rac$ is finite, the above sequence terminates for a fixed $w$.\\
\noindent \textbf{Remark 2} $n$ can be found straightforwardly from the $\hbar$ series of $F_{C_{(r,2)}(\bm{4_1})}(x,q)$. We propose a formula for $n$ below.\\
\noindent \textbf{Remark 3} A power of each $q$-monomial increases by multiples for a positive integer from $f^{\pm}_{m}$ to  $f^{\pm}_{m+4}$.\\
\noindent \textbf{Remark 4} Coefficient functions $\lac f_d (q) \rac$ that do not fit into the above pattern vanish.\\
\noindent \textbf{Remark 5} $h_1 (q) q^{s}$ is added when $h_{4w+3}$ appears in $f^{\pm }_{b}(q)$ for any $b$.\\
\newline


\begin{conjecture} 
For the above family of knots, the power $n \in \intg$ in Conjecture 2.1 is given by $w+1$. 
\end{conjecture}

Since the above $f_m(q)$ are polynomials rather than series, we have the following consequence for the mirror knot $m(K)$.
\begin{cor} 
For a $(2w+1,2)$-cabling of $\bm{4_1}$ where $w < -3$, its nonzero initial data $\lac f_m (q) \rac = \lac f^{+}_i (q) \rac \bigsqcup \lac f^{-}_j (q) \rac$ for a recursion relation are given by $q \longleftrightarrow q^{-1}$ in Conjecture 2.1.
\end{cor}

We state that the patterns in Conjecture 2.1 holds for all nonzero $f_m (q)$.

\begin{conjecture} 
For the family of the cable knots $K_w =C_{(2w+1,2)}(\bm{4_1})$ in $\intg HS^3$ for $w>3$, all of their coefficient functions $f^{\pm}_m (q) \in \intg[q^{\pm 1}]$ of $F_{K_w}(x,q)$ can be expressed as linear combinations of the coefficient functions $\lac h_s (q) \rac$ of $F_{\bm{4_1}}(x,q)$ following the pattern in Conjecture 2.1.
\end{conjecture}

We next verify Conjecture 2.1 and 2.4 for some examples. 

\section{Examples}

\subsection{Example 1} We first analyze $K=C_{(9,2)}(\bm{4_1})$.  Applying Conjecture 2.1 to $w=4$, its predictions agree with the initial data listed in Section 7 of \cite{C}\footnote{There is an overall factor of $2$ difference due to a different convention of $f_m (q)$ in $F_K$}. 
\newline

\subsection{Example 2} The next simplest case is $K=C_{(11,2)}(\bm{4_1})$.
\begin{prop}
The $\hbar$ expansion of the $J_{K,n}(q)$  is given by
\begin{equation*}
\begin{aligned}		
	J_{K,n}(e^{\hbar}) & = 1 + \lb 11 - 11 n^2 \rb \hbar^2 + \lb -88 + 88 n^2 \rb \hbar^3 + \lb \frac{11891}{12} - 1137 n^2 + \frac{1753}{12} n^4 \rb \hbar^4\\
	& + \lb - 12826 + \frac{47036}{3} n^2 - \frac{8558}{3} n^4 \rb \hbar^5 + \lb  \frac{69672971}{360} - \frac{991683}{4} n^2 + \frac{224545}{4} n^4 - \frac{630551}{360} n^6 \rb \hbar^6 + \cdots
\end{aligned}
\end{equation*}
\end{prop}
\noindent At each $\hbar$ order, the degree of the polynomial in $n$ is at most the order of $\hbar$, which is an equivalent characterization of the MMR expansion of the colored Jones polynomial of a knot.
\newline

The cabling formula for the Alexander polynomial of a knot $K$ is~\cite{Hedden}
$$
\Delta_{C_{(p,q)}(K)}(t)= \Delta_K (t^p) \Delta_{T_{(p,q)}}(t), \quad 2 \leq p < |q| \quad \text{gcd}(p,q)=1, 
$$
where $\Delta(t)$ is the symmetrized Alexander polynomial and $T_{(p,q)}$ is the $(p,q)$ torus knot. Note that our convention for the parameters of the torus knot are switched (i.e. $p\equiv 2, q\equiv r$). 
\begin{lemma}
The Alexander polynomial of $C_{(11,2)}(\bm{4_1})$ is as follows.
\begin{equation*}
\begin{aligned}		
\Delta_{C_{(11,2)}{(\bm{4_1})}}(x) & = \Delta_{\bm{4_1}} (x^2) \Delta_{T_{(2,11)}}(x)\\
& = -1-\frac{1}{x^7}+\frac{1}{x^6}+\frac{2}{x^5}-\frac{2}{x^4}+\frac{1}{x^3}-\frac{1}{x^2}+\frac{1}{x}+x-x^2+x^3-2 x^4+2 x^5+x^6-x^7.
\end{aligned}
\end{equation*}
\end{lemma}
\noindent From this Alexander polynomial its symmetric expansion about $x=0$ (in $x$) and $x=\infty$ (in $1/x$) can be computed.
\begin{prop} 
The symmetric expansion of the inverse of the Alexander polynomial of $K$ in the limit of $\hbar \rarw 0$ is given by
\begin{align*}
\lim_{q \rarw 1} 2F_{K}(x,q) &  =  2\, \text{s.e.} \lb \frac{x^{1/2}-x^{-1/2}}{\Delta_{K}(x)} \rb\\
 & =  x^{13/2}-\left(\frac{1}{x}\right)^{13/2}+2 x^{17/2}-2 \left(\frac{1}{x}\right)^{17/2}+5x^{21/2}-5 \left(\frac{1}{x}\right)^{21/2}+13 x^{25/2}-13 \left(\frac{1}{x}\right)^{25/2}\\
 & +34   x^{29/2}-34 \left(\frac{1}{x}\right)^{29/2}+89 x^{33/2}-89\left(\frac{1}{x}\right)^{33/2}-x^{35/2}+\left(\frac{1}{x}\right)^{35/2}+233 x^{37/2}-233\left(\frac{1}{x}\right)^{37/2}\\
 & -2 x^{39/2}+2 \left(\frac{1}{x}\right)^{39/2}+610 x^{41/2}-610\left(\frac{1}{x}\right)^{41/2}-5 x^{43/2}+5 \left(\frac{1}{x}\right)^{43/2} + \cdots \in \intg \lsb \lsb x^{\pm 1/2} \rsb \rsb 
\end{align*}
where \text{s.e.} denotes the symmetric expansion in $x$ about $x=0$ and $1/x$ about $x=\infty$.
\end{prop} 
\noindent The coefficients in the expansions are integers and the Alexander polynomial is monic, which is a necessary condition for $f_m(q)$'s in (1) to be polynomials. Furthermore, the coefficients of $x^{35/2}, x^{39/2}, x^{43/2}$ etc are negative, which are reflected as negative coefficient elements in the initial data below (cf.Proposition 3.5).
\newline


\noindent From (5) we find the recursion relation for $f_{m}$.
\begin{theorem}
The recursion relation for $f_{m}(q) \in \intg[q^{\pm 1}]$ of $F_K (x,q)$ is given by
\begin{equation}
\begin{aligned}
f_{v+102}\,(q) & = \frac{-1}{q^{\frac{115+v}{2}} \left(1-q^{\frac{89+v}{2}}\right)} \Big[ t_2 f_{98+v}  + t_4 f_{94+v}  + t_6 f_{90+v} + t_8 f_{86+v}  + t_{10} f_{82+v} +t_{11} f_{80+v} \\
& + t_{12} f_{78+v}   + t_{13} f_{76+v}  + t_{14} f_{74+v} + t_{15} f_{72+v}  + t_{16} f_{70+v} +  t_{17} f_{68+v} + t_{18} f_{66+v} \\
& +t_{19} f_{64+v}  + t_{20} f_{62+v}  +  t_{21} f_{60+v} +  t_{22} f_{58+v} +  t_{23} f_{56+v} + t_{24} f_{54+v}  + t_{25} f_{52+v} \\
&  + t_{26} f_{50+v}  + t_{27} f_{48+v}  + t_{28} f_{46+v}  + t_{29} f_{44+v}  + t_{30}  f_{42+v}+  t_{31} f_{40+v} + t_{32} f_{38+v}\\
& + t_{33}  f_{36+v} + t_{34} f_{34+v}  +  t_{35} f_{32+v} + t_{36} f_{30+v}  + t_{37} f_{28+v}  + t_{38} f_{26+v}  + t_{39} f_{24+v} \\
&  +  t_{40} f_{22+v} + t_{41} f_{20+v}  + t_{43} f_{16+v}  +  t_{45} f_{12+v} +  t_{47} f_{8+v} +  t_{49} f_{4+v} + t_{51} f_v \vphantom{1} \Big] 
\end{aligned}
\end{equation}
\noindent where $t_{s}=t_{s}(q,q^v)$ are $\intg$-Laurent polynomials in $q$ and $q^v$. 
\end{theorem}

\begin{prop}
The initial data for (6) consists of $\lac f_{j} (q) \in \intg[q^{\pm 1}] | j=1,\cdots ,101 \rac$. All the nonzero elements are the following.
\begin{align*}	
f^{+}_{13} (q) & = 2\,  h_1 (q)\, q^6 \\
f^{+}_{17} (q) & = 2\, h_3 (q)\, q^{7} \\
f^{+}_{21} (q) & = 2\, h_5 (q)\, q^{8}\\
&\vdots\\
f^{+}_{57} (q) & = 2\, \lb  h_{23} (q)\, q^{17 } + h_1 (q) \, q^{41}  \rb\\
f^{+}_{61} (q) & = 2\, \lb h_{ 25 } (q)\, q^{18} + h_3 (q) \, q^{46 }  \rb\\ 
f^{+}_{65} (q) & = 2\, \lb h_{27 } (q)\, q^{19} + h_5 (q) \, q^{51 }  \rb\\
&\vdots\\
f^{+}_{101} (q) & = 2\, \lb h_{ 45} (q)\, q^{28} + h_{23 } (q) \, q^{96 }  + h_1 (q) \, q^{120 }  \rb\\
&\vdots\\
f^{-}_{35} (q) & = -2\,  h_1 (q)\, q^{18} \\
f^{-}_{39} (q) & = -2\, h_3 (q)\, q^{21}  \\
f^{-}_{43} (q) & = -2\, h_5 (q)\, q^{24}\\
&\vdots\\
f^{-}_{79} (q) & = -2 \lb h_{23} \, q^{51} + h_1 \, q^{75} \rb\\
f^{-}_{83} (q) & = -2 \lb h_{25} \, q^{54} + h_3 (q) \, q^{82} \rb\\
&\vdots\\
f^{-}_{99} (q) & = -2 \lb h_{33} \, q^{66} + h_{11} (q) \, q^{110} \rb\\
\end{align*}
\end{prop}
After substituting them into (6), we obtain $f_{103} (q)$ and $f_{105} (q)$ as Laurent polynomials with integer coefficients, which is non-trivial because of the denominator $(1-q^{(89+v)/2})$ in (6). Recasting them in terms of $h_s (q)$, we get
\begin{align}
f^{-}_{103} (q) & = -2 \lb h_{35}(q) \, q^{69} + h_{13} (q) \, q^{117} \rb\\
f^{+}_{105} (q) & = 2\, \lb h_{ 47} (q)\, q^{29} + h_{25 } (q) \, q^{101}  + h_3 (q) \, q^{129}  \rb
\end{align}
They agree with the predictions from Conjecture 2.4 by extrapolating it beyond $f_{101}(q)$. As another check, $\hbar$ expansions of (7) and (8) match with that of $F_{C_{(11,2)}(\bm{4_1})}(x,e^{\hbar})$. They can in turn be used to confirm that subsequent $f_m (q)$'s obtained from Conjecture 2.4 agree with that of (6).
\begin{align*}
f^{-}_{107} (q) & = -2 \lb h_{37}(q) \, q^{72} + h_{15} (q) \, q^{124} \rb\\
f^{+}_{109} (q) & = 2\, \lb h_{ 49} (q)\, q^{30} + h_{27 } (q) \, q^{106}  + h_5 (q) \, q^{138}  \rb\\
\end{align*}
Furthermore, their $\hbar$ expansions coincide with that of $F_{C_{(13,2)}(\bm{4_1})}$.
\newline

\noindent \textbf{Remark} When $h_{23}(q)$ appears in $f^{\pm}_b$ for any $b$, $h_1(q)$ is added.
\newline

\subsubsection{Dehn surgery and $\hat{Z}$}

In \cite{GPPV}, it was conjectured that the $sl(2)$ WRT invariant of a rational homology 3-sphere $Y\, (b_1(Y)=0)$ decomposes into a linear combination of power series $\hat{Z}(q)$, which converges in $|q|<1$ and itself is an invariant of $Y$. Specifically,
$$
WRT[Y;k]= \frac{1}{i\sqrt{2k}} \sum_{a,b \in Spin^c (Y)} c^{WRT}_{ab}\, \hat{Z}_{b}(q) \Big\vert_{q \rarw e^{\frac{i2\pi}{k}}},
$$
where
$$
\hat{Z}_{b}(Y;q) \in \frac{1}{2^c}q^{\Delta_{b}} \intg[[q]], \qquad |q|<1,\qquad \Delta_{b} \in \rational,\quad c \in \intg_{\geq 0}.
$$
Performing a Dehn surgery on a complement of a knot $K$ along a (negative) slope $-p/r$~\footnote{Surgery along a positive slope may not produce a power series; there has been some progress on positive surgery in \cite{P}.} yields a relation between $F_K (x,q)$ and $\hat{Z}(q)$ given by~\cite{GM}
$$
\hat{Z}_{b}[Y_{-\frac{p}{r}}(K);q] = \pm q^d \mathcal{L}^{(b)}_{-p/r} \lsb \lb x^{\frac{1}{2r}} - x^{-\frac{1}{2r}} \rb F_{K}(x,q) \rsb \qquad d \in \rational,
$$
$$
\mathcal{L}^{(b)}_{-p/r} : x^{u}q^{v} \mapsto \begin{cases}
      q^{u^2 r/p}q^v & \text{if}\quad ru - b \in p\intg \\
      0 & \text{otherwise}
\end{cases}
$$
Performing surgeries along $-1/r$ slopes for a few $r$ results in
\begin{align*}
\hat{Z}_0 [S^{3}_{-\frac{1}{2}}(C_{11,2}(4_1));q] & = q^{\frac{167}{2}} \Big( -1+q^{13}-2 q^{59}+2 q^{76}-q^{133}-3 q^{134}-q^{135}+q^{154}+3 q^{155}+q^{156}\\
& -2q^{223}-2 q^{224}-5 q^{225}-2 q^{226}-2 q^{227}+2 q^{248}+2 q^{249} + \cdots  \Big)\\
\hat{Z}_0 [S^{3}_{-\frac{1}{3}}(C_{11,2}(4_1));q] & = q^{\frac{251}{2}}  \Big( -2+2 q^{13}-4 q^{89}+4 q^{106}-2 q^{201}-6 q^{202}-2 q^{203}+2 q^{222}+6 q^{223}\\
& +2 q^{224} -4 q^{337}-4 q^{338}-10 q^{339}-4 q^{340}-4 q^{341}+4 q^{362}+4 q^{363} + \cdots  \Big)\\
\hat{Z}_0 [S^{3}_{-\frac{1}{4}}(C_{11,2}(4_1));q] & = q^{\frac{335}{2}} \Big( -2+2 q^{13}-4 q^{119}+4 q^{136}-2 q^{269}-6 q^{270}-2 q^{271}+2 q^{290}+6 q^{291}\\
& +2 q^{292} -4q^{451}-4 q^{452}-10 q^{453}-4 q^{454}-4 q^{455}+4 q^{476}+4 q^{477} + \cdots  \Big)\\
\end{align*}
\newline

\subsection{Example 3 } We next consider $K=C_{(13,2)}(\bm{4_1})$.

\begin{prop}
The $\hbar$ expansion of the $J_{K,n}(q)$  is given by
\begin{equation*}
\begin{aligned}		
	J_{K,n}(e^{\hbar}) & = 1 + \lb 17 - 17 n^2 \rb \hbar^2 + \lb -156 + 156 n^2 \rb \hbar^3 + \lb \frac{24749}{12} - 2365 n^2 + \frac{3631}{12} n^4 \rb \hbar^4\\
	& + \lb - 31629 + 38662 n^2 - 7033 n^4 \rb \hbar^5 + \Big(  \frac{203413517}{360} - \frac{8687953}{12} n^2 + \frac{1969367}{12} n^4 \\
	& - \frac{1855937}{360} n^6 \Big) \hbar^6 + \cdots
\end{aligned}
\end{equation*}
\end{prop}
\noindent This confirms that the MMR condition is satisfied by $K$. 

\begin{lemma}
The Alexander polynomial of $C_{(13,2)}(\bm{4_1})$ is given by
\begin{equation*}
\begin{aligned}		
\Delta_{C_{(13,2)}{(\bm{4_1})}}(x) & = \Delta_{\bm{4_1}} (x^2) \Delta_{T_{(2,13)}}(x)\\
& = -t^8-\frac{1}{t^8}+t^7+\frac{1}{t^7}+2 t^6+\frac{2}{t^6}-2t^5-\frac{2}{t^5}+t^4+\frac{1}{t^4}-t^3-\frac{1}{t^3}+t^2+\frac{1}{t^2}-t-\frac{1}{t}+1.
\end{aligned}
\end{equation*}
\end{lemma}

\begin{prop}
The symmetric expansion of the inverse of the Alexander polynomial of $K$ in the limit of $\hbar \rarw 0$ is given by
\begin{align*}
\lim_{q \rarw 1} 2F_{K}(x,q) &  =  2\, \text{s.e.} \lb \frac{x^{1/2}-x^{-1/2}}{\Delta_{K}(x)} \rb\\
 & = x^{15/2}-\left(\frac{1}{x}\right)^{15/2}+2 x^{19/2}-2 \left(\frac{1}{x}\right)^{19/2}+5 x^{23/2}-5\left(\frac{1}{x}\right)^{23/2}+13 x^{27/2}-13 \left(\frac{1}{x}\right)^{27/2}\\
 & +34 x^{31/2}-34\left(\frac{1}{x}\right)^{31/2}+89 x^{35/2}-89 \left(\frac{1}{x}\right)^{35/2}+233 x^{39/2}-233\left(\frac{1}{x}\right)^{39/2}-x^{41/2}+\left(\frac{1}{x}\right)^{41/2}\\
 & +610 x^{43/2}-610\left(\frac{1}{x}\right)^{43/2}-2 x^{45/2}+2 \left(\frac{1}{x}\right)^{45/2} + \cdots \in \intg \lsb \lsb x^{\pm 1/2} \rsb \rsb 
\end{align*}
\end{prop} 
\noindent The coefficients of $x^{41/2}$ and $x^{45/2}$ etc being negative is consistent with the negative coefficient elements in the initial data below (cf. Proposition 3.10).
\newline

\begin{theorem}
The recursion relation for $f_{m}(q) \in \intg[q^{\pm 1}]$ of the above $F_K (x,q)$ is given by
\begin{equation}
\begin{aligned}
f_{v+106}\,(q) & = \frac{-1}{q^{\frac{121+v}{2}} \left(1-q^{\frac{91+v}{2}}\right)} \Big[ t_2 f_{v+102}+t_4 f_{v+98}+t_6 f_{v+94}+t_8 f_{v+90}+t_{10} f_{v+86}+ t_{12} f_{v+82}\\
& +t_{13} f_{v+80}+t_{14} f_{v+78}+t_{15} f_{v+76}+t_{16} f_{v+74}+t_{17} f_{v+72}+t_{18} f_{v+70}+t_{19} f_{v+68}\\
& +t_{20} f_{v+66}+t_{21} f_{v+64}+t_{22} f_{v+62}+t_{23} f_{v+60}+t_{24} f_{v+58}+t_{25}  f_{v+56}+t_{26} f_{v+54}\\
& +t_{27} f_{v+52}+t_{28} f_{v+50}+t_{29} f_{v+48}+t_{30} f_{v+46}+t_{31} f_{v+44}+t_{32} f_{v+42}+t_{33} f_{v+40} \\
& +t_{34} f_{v+38}+t_{35} f_{v+36}+t_{36} f_{v+34}+t_{37} f_{v+32}+t_{38} f_{v+30}+t_{39} f_{v+28}+t_{40} f_{v+26}\\
& +t_{41} f_{v+24}+t_{43} f_{v+20}+t_{45}f_{v+16}+t_{47} f_{v+12}+t_{49} f_{v+8}+t_{51} f_{v+4}+t_{53} f_v \vphantom{1} \Big] 
\end{aligned}
\end{equation}
where $t_{s}=t_{s}(q,q^v)$ are $\intg$-Laurent polynomials in $q$ and $q^v$.
\end{theorem}

\begin{prop}
The initial data for (9) consists of $\lac f_{j} (q) \in \intg[q^{\pm 1}] | j=1,\cdots ,105 \rac$. All the nonzero elements of the set are the following.
\begin{align*}	
f^{+}_{15} (q) & = 2\,  h_1 (q)\, q^7 \\
f^{+}_{19} (q) & = 2\, h_3 (q)\, q^8  \\
f^{+}_{23} (q) & = 2\, h_5 (q)\, q^9\\
&\vdots\\
f^{+}_{67} (q) & = 2\, \lb  h_{27} (q)\, q^{20 } + h_1 (q) \, q^{48}  \rb\\
f^{+}_{71} (q) & = 2\, \lb h_{ 29 } (q)\, q^{21} + h_3 (q) \, q^{53}  \rb\\ 
&\vdots\\
f^{-}_{41} (q) & = -2\,  h_1 (q)\, q^{21} \\
f^{-}_{45} (q) & = -2\, h_3 (q)\, q^{24}  \\
f^{-}_{49} (q) & = -2\, h_5 (q)\, q^{27}\\
&\vdots\\
f^{-}_{93} (q) & = -2 \lb h_{27}(q) \, q^{60} + h_1 \, q^{88} \rb\\
f^{-}_{97} (q) & = -2 \lb h_{29}(q) \, q^{63} + h_3 (q) \, q^{95} \rb\\
&\vdots\\
f^{-}_{105} (q) & = -2 \lb h_{33}(q) \, q^{69} + h_{7} (q) \, q^{109} \rb\\
\end{align*}
\end{prop}
\noindent After substituting them into (9), we indeed obtain $f_{107} (q)$ and $f_{109} (q)$ as Laurent polynomials with $\intg$-coefficients. They can be recast in terms of $h_s (q)$:
\begin{align}
f^{+}_{107} (q) & = 2 \lb h_{47}(q) \, q^{30} + h_{21} (q) \, q^{98} \rb\\
f^{-}_{109} (q) & = -2\, \lb h_{ 35} (q)\, q^{72} + h_{9 } (q) \, q^{116}   \rb
\end{align}
The predictions from Conjecture 2.4 coincide with the above. Furthermore, $\hbar$ expansions of (10) and (11) agree with that of $F_{C_{(13,2)}(\bm{4_1})}(x,e^{\hbar})$. These can in turn be used to confirm that subsequent $f_m (q)$'s obtained from (9) again agree with that of Conjecture 2.4.
\begin{align*}
f^{+}_{111} (q) & = 2 \lb h_{49}(q) \, q^{31} + h_{23} (q) \, q^{103} \rb\\
f^{-}_{113} (q) & = -2\, \lb h_{ 37} (q)\, q^{75} + h_{11 } (q) \, q^{123}   \rb\\
\end{align*}
Their $\hbar$ expansions of the above results match with that of $F_{C_{(13,2)}(\bm{4_1})}$ as well.
\newline

\noindent \textbf{Remark} When $h_{27}(q)$ appears in $f^{\pm}_b$ for any $b$, $h_1(q)$ is added.
\newline

\section{$(3w+1,3)$-cables\, $(w>3)$}

The quantum (or noncommutative) A-polynomial of a class of cable knot $C_{(r,3)}(\bm{4_1})$ is qualitatively different from that of the $(2w+1,2)$-cabling in Section 2. It is given by~\cite{Ru}
\begin{equation}
\hat{A}_{K} (t,M,L)=(L-1)B(t,M)^{-1}Q(t,M,L) \lb t^{6s} M^{3r} L^2 - t^{-6r}M^{-3r}, \rb \quad gcd(r,3)=1, \quad |r| > 12.
\end{equation}
The definitions of $Q(t,M,L)$ and $B(t,M)$ are written in Appendix A. Using (3) together with (12) yields
\begin{equation}
\begin{split}
\phantom{A} & \phantom{=} c_0 (x,q)F_K(x,q) + c_1 (x,q)F_K(xq,q)+ c_2 (x,q)F_K(xq^2,q) + c_3 (x,q) F_K(xq^3,q)\\
 & + c_4 (x,q) F_K(xq^4,q) + c_5 (x,q) F_K(xq^5,q)=0,
\end{split}
\end{equation} 
where $c_0  ,c_1 ,c_2 , c_3 , c_4 , c_5  $ are rational functions. We state the nonzero initial data needed for (13).

\begin{conjecture} 
For a $(3w+1,3)$-cabling of $\bm{4_1}\, (w > 3)$ in $\intg HS^3$, nonzero initial data $\lac f_m (q) \in \intg [q^{\pm 1}] \rac = \lac f^{+}_i (q)  |\, \text{all coefficients are positive} \rac \bigsqcup \lac f^{-}_j (q) |\, \text{all coefficients are negative}  \rac$ for (13) are determined by the coefficient functions $\lac h_s (q) \rac$ of  $F_{\bm{4_1}}(x,q)$ as follows.
\begin{align*}	
f^{+}_{6w+5} (q) & = 2\,  h_1 (q)\, q^n \qquad ( h_1 (q)=1)\\
f^{+}_{6w+11} (q) & = 2\, h_3 (q)\, q^{n+2}  \qquad ( h_3 (q)=2)\\
f^{+}_{6w+17} (q) & = 2\, h_5 (q)\, q^{n+4}\\
&\vdots\\
f^{+}_{24w+11} (q) & = 2\, \lb  h_{6w+3 } (q)\, q^{n+ 6 w + 2  } + h_1 (q) \, q^{n+ 6 w + 2 +  \Delta}  \rb\\
f^{+}_{24w+17} (q) & = 2\, \lb h_{ 6w+5 } (q)\, q^{n + 6 w + 4 } + h_3 (q) \, q^{n + 6 w + 10 + \Delta }  \rb\\ 
f^{+}_{24w+23} (q) & = 2\, \lb h_{ 6w+7 } (q)\, q^{n+  6 w + 6 } + h_5 (q) \, q^{n + 6 w + 18  + \Delta }  \rb\\
&\vdots\\
f^{+}_{42w+17} (q) & = 2\, \lb h_{ 12w+5  } (q)\, q^{n+  12w + 4 } + h_{6w+3 } (q) \, q^{n + 30w +10 +  \Delta }  + h_1 (q) \, q^{n +  30w +10 + 2\Delta }  \rb\\
f^{+}_{42w+23} (q) & = 2\, \lb h_{12w+7} (q)\, q^{n+12w +6 } + h_{6w+5} (q) \, q^{n + 30w +18 +  \Delta }  + h_3 (q)\, q^{n + 30w +24 + 2\Delta }  \rb\\
&\vdots\\
f^{+}_{60w+23} (q) & = 2\, \lb h_{ 18w+ 7} (q)\, q^{n+ 18w + 6} + h_{12w+5 } (q) \, q^{n + 54w +18  + \Delta }  + h_{6w+3} (q) \, q^{n + 72w+24  + 2\Delta } \notag\right. \\
& \notag\left. + h_{1} (q) \, q^{n + 72w+24  + 3\Delta } \rb\\
f^{+}_{60w+29} (q) & = 2\, \lb h_{18w+9} (q)\, q^{n+12w +8 } + h_{12w+7} (q) \, q^{n + 54w +26 +  \Delta } + h_{6w+5} (q) \, q^{n + 72w + 38  + 2\Delta }\notag\right.\\
& \notag\left. + h_3 (q)\, q^{n + 72w + 44 + 3\Delta }  \phantom{1} \rb\\
&\vdots\\
f^{+}_{78w+29} (q) & = 2\, \lb h_{24w+9} (q)\, q^{n+ 24w + 8} + h_{18w+7 } (q) \, q^{n + 78w + 26 + \Delta }  + h_{12w+5} (q) \, q^{n + 114w + 38 + 2\Delta} \notag\right. \\
& \notag\left. + h_{6w+3} (q) \, q^{n + 132w + 44 + 3\Delta } + h_{1} (q) \, q^{n + 132w + 44 + 4\Delta }  \phantom{1}\rb\\
f^{+}_{78w+35} (q) & = 2\, \lb h_{24w+11} (q)\, q^{n+ 24w + 10} + h_{18w+9 } (q) \, q^{n + 78w + 34 + \Delta }  + h_{12w+7} (q) \, q^{n + 114w + 52 + 2\Delta } \notag\right.\\
& \notag\left. + h_{6w+5} (q) \, q^{n + 132w + 64 + 3\Delta } + h_{3} (q) \, q^{n + 132w + 70 + 4\Delta }  \phantom{1}\rb\\
&\vdots\\
f^{-}_{12w+7} (q) & = -2\,  h_1 (q)\, q^{2n} \\
f^{-}_{12w+13} (q) & = -2\, h_3 (q)\, q^{2n+4}  \\
f^{-}_{12w+19} (q) & = -2\, h_5 (q)\, q^{2n+8}\\
&\vdots\\
f^{-}_{30w+13} (q) & = -2 \lb h_{6w+3}(q) \, q^{2n+12 w + 4  } + h_1 (q) \, q^{2n+ 12 w + 4   + \Delta} \rb\\
f^{-}_{30w+19} (q) & = -2 \lb h_{6w+5}(q) \, q^{2n+12 w + 8 } + h_3 (q) \, q^{2n+ 12 w + 14  + \Delta} \rb\\
&\vdots\\
f^{-}_{48w+19} (q) & = -2\, \lb h_{12w+5} (q)\, q^{2n+ 24w +8 } + h_{6w+3 } (q) \, q^{2n + 42w + 14 +  \Delta }  + h_1 (q) \, q^{2n + 42w + 14 + 2\Delta }  \rb\\
f^{-}_{48w+25} (q) & = -2\, \lb h_{12w+7} (q)\, q^{2n+ 24w +12 } + h_{6w+5 } (q) \, q^{2n + 42w + 24 +  \Delta }  + h_3 (q) \, q^{2n +  42w +30 + 2\Delta }  \rb\\
&\vdots\\
f^{-}_{66w+25} (q) & = -2\, \lb h_{18w+7} (q)\, q^{2n+ 36w +12 } + h_{12w+5 } (q) \, q^{2n + 72w + 24 +  \Delta }  + h_{6w+3} (q) \, q^{2n + 90w +30 + 2\Delta } \notag\right.\\
& \notag\left. + h_1 (q) \, q^{2n + 90w +30 + 3\Delta }  \phantom{1}\rb\\
f^{-}_{66w+31} (q) & = -2\, \lb h_{18w+9} (q)\, q^{2n+ 36w +16 } + h_{12w+7 } (q) \, q^{2n + 72w + 34 +  \Delta }  + h_{6w+5} (q) \, q^{2n + 90w + 46 + 2\Delta } \notag\right.\\
& \notag\left. + h_3 (q) \, q^{2n + 90w +52 + 3\Delta }  \phantom{1}\rb\\
&\vdots\\
\end{align*}
where $n \in \intg$ and $\Delta =9w+6$.
\end{conjecture}
\noindent We state several remarks.\\
\noindent \textbf{Remark 1} As in the previous cabling, the above sequence terminates for a fixed $w$.\\
\noindent \textbf{Remark 2} $n$ can be found easily from the $\hbar$ series of $F_{C_{(r,2)}(\bm{4_1})}(x,q)$. We propose a formula for $n$ below.\\
\noindent \textbf{Remark 3} A power of each $q$-monomial increases by multiples for a positive integer from $f^{\pm}_{m}$ to  $f^{\pm}_{m+6}$\\
\noindent \textbf{Remark 4} Coefficient functions $ f_d (q) $ that do not fit into the above pattern vanish.\\
\noindent \textbf{Remark 5} $h_1 (q) q^{s}$ is added when $h_{6w+3}$ appears in $f^{\pm}_{b}(q)$ for any $b$. 
\newline

\begin{conjecture} 
For the above family of knots, the power $n \in \intg$ in Conjecture 4.1 is given by $3w+2$. 
\end{conjecture}
\noindent As in the $(2w+1,2)$-cabling case, we get $f_m(q)$ for the mirror knot $m(K)$.
\begin{cor} 
For a $(3w+1,3)$-cabling of $\bm{4_1}$ where $w < -3$, its nonzero initial data $\lac f_m (q) \rac = \lac f^{+}_m (q) \rac \bigsqcup \lac f^{-}_m (q) \rac$ for a recursion relation are given by $q \longleftrightarrow q^{-1}$ in Conjecture 4.1.
\end{cor}

\begin{conjecture} 
For the family of the cable knots $K_w =C_{(3w+1,3)}(\bm{4_1})$ in $\intg HS^3$ for $w>3$, all of their coefficient functions $f^{\pm}_m (q) \in \intg[q^{\pm 1}]$ of $F_{K_w}(x,q)$ can be expressed as linear combinations of the coefficient functions $\lac h_s (q) \rac$ of $F_{\bm{4_1}}(x,q)$ following the pattern of Conjecture 4.1.
\end{conjecture}

\section{Examples}

\subsection{Example 1} We begin with the simplest example $K=C_{(13,3)}(\bm{4_1})$.

\begin{prop}
The $\hbar$ expansion of the $J_{K,n}(q)$  is given by
\begin{align*}
	J_{K,n}(e^{\hbar}) & = 1 + \lb 47 - 47 n^2 \rb \hbar^2 + \lb -624 + 624 n^2 \rb \hbar^3 + \lb \frac{151919}{12} - 14605 n^2 + \frac{23341}{12} n^4 \rb \hbar^4\\
	& + \lb - 294528 + 361088 n^2 - 66560 n^4 \rb \hbar^5 + \Big(  \frac{2864712407}{360} - \frac{122607733}{12} n^2 + \frac{28027787}{12} n^4 \\
	& - \frac{27314027}{360} n^6 \Big) \hbar^6 + \cdots
\end{align*}
\end{prop}

\begin{lemma}
The Alexander polynomial of $C_{(13,3)}(\bm{4_1})$ is given by
\begin{equation*}
\begin{aligned}		
\Delta_{C_{(13,3)}{(\bm{4_1})}}(x) & = \Delta_{\bm{4_1}} (x^3) \Delta_{T_{(3,13)}}(x)\\
& = -t^{15}-\frac{1}{t^{15}}+t^{14}+\frac{1}{t^{14}}+2 t^{12}+\frac{2}{t^{12}}-2 t^{11}-\frac{2}{t^{11}}+t^9+\frac{1}{t^9}-t^8-\frac{1}{t^8}+t^6+\frac{1}{t^6}\\
& -t^5-\frac{1}{t^5}+t^3+\frac{1}{t^3}-2 t^2-\frac{2}{t^2}+t+\frac{1}{t}+1\\
\end{aligned}
\end{equation*}
\end{lemma}

\begin{prop}
The symmetric expansion of the inverse of the Alexander polynomial of $K$ in the limit of $\hbar \rarw 0$ is given by
\begin{align*}
\lim_{q \rarw 1} 2F_{K}(x,q) &  =  2\, \text{s.e.} \lb \frac{x^{1/2}-x^{-1/2}}{\Delta_{K}(x)} \rb\\
& = 2 x^{29/2} -\frac{2}{x^{29/2}} +4 x^{35/2} -\frac{4}{x^{35/2}} +10 x^{41/2} -\frac{10}{x^{41/2}} +26 x^{47/2} -\frac{26}{x^{47/2}} +68 x^{53/2}\\
& -\frac{68}{x^{53/2}} -2 x^{55/2} +\frac{2}{x^{55/2}} +178 x^{59/2} -\frac{178}{x^{59/2}} -4 x^{61/2} +\frac{4}{x^{61/2}} +466x^{65/2} \\
& -\frac{466}{x^{65/2}}-10 x^{67/2}+\frac{10}{x^{67/2}} +1220 x^{71/2}-\frac{1220}{x^{71/2}} + \cdots \in \intg \lsb \lsb x^{\pm 1/2} \rsb \rsb 
\end{align*}
\end{prop}
\noindent In this case, the coefficients of $x^{55/2}, x^{61/2}, x^{67/2}$ and so forth are negative, which are reflected as negative coefficient elements in the initial data.
\newline

\begin{prop}
The initial data for the recursion relation of $K$ consists of $I=\lac f_{j} (q) \in \intg[q^{\pm 1}] | j=1,\cdots ,353 \rac$. The nonzero elements in $I$ are obtained by setting $w=4$ and hence $n=14$ in Conjecture 4.1.
\end{prop}

\noindent We extrapolate Conjecture 4.1 beyond $I$.
\begin{align*}
f^{-}_{355}(q) & = -2 \lb h_{101}\, q^{228} + h_{75}\, q^{492} + h_{49}\, q^{678} + h_{23}\, q^{786} \rb\\
f_{357}(q) & = 0
\end{align*} 
$f^{-}_{355}$ indeed agree with that of the recursion result and the vanishing of $f_{357}$ is consistent with the recursion method as well. We confirm that several subsequent $f_m (q)$'s match with the results of the recursion. For example,
\begin{align*}
f^{+}_{359}(q) & = 2 \lb h_{111}\, q^{124} + h_{85}\, q^{418}+ h_{59}\, q^{634}+ h_{33}\, q^{772} + h_7 \, q^{832} \rb \\
f^{-}_{361}(q) & = -2 \lb h_{103}\, q^{232} +  h_{77}\, q^{502} + h_{51}\, q^{694} + h_{25}\, q^{808} \rb 
\end{align*} 
\newline

We perform Dehn surgeries along $-1/r$ slopes for a few $r$.
\begin{align*}
\hat{Z}_0 [S^{3}_{-1}(C_{13,3}(4_1));q] & = q^{\frac{419}{2}} \Big( -2+2 q^{29}-4 q^{95}+4 q^{130}-2 q^{207}-6 q^{208}-2 q^{209}+2 q^{248}+6 q^{249}\\
& +2   q^{250}-4 q^{337}-4 q^{338}-10 q^{339}-4 q^{340}-4 q^{341}+4 q^{384}+4 q^{385} + \cdots  \Big)\\
\hat{Z}_0 [S^{3}_{-\frac{1}{2}}(C_{13,3}(4_1));q] & = q^{\frac{839}{2}}  \Big( -2+2 q^{29}-4 q^{191}+4 q^{226}-2 q^{417}-6 q^{418}-2 q^{419}+2 q^{458}\\
& +6 q^{459}+2 q^{460}-4 q^{679}-4 q^{680}-10 q^{681}-4 q^{682}-4 q^{683}+4 q^{726}+4 q^{727} + \cdots  \Big)\\
\hat{Z}_0 [S^{3}_{-\frac{1}{3}}(C_{13,3}(4_1));q] & = q^{\frac{1259}{2}} \Big( -2+2 q^{29}-4 q^{287}+4 q^{322}-2 q^{627}-6 q^{628}-2 q^{629}+2 q^{668}\\
& +6 q^{669}+2 q^{670}-4 q^{1021}-4 q^{1022}-10 q^{1023}-4 q^{1024}-4 q^{1025}+4 q^{1068}+ \cdots  \Big)\\
\end{align*}
\newline

\subsection{Example 2} We move onto $K=C_{(16,3)}(\bm{4_1})$.	

\begin{prop}
The $\hbar$ expansion of the $J_{K,n}(q)$  is given by
\begin{align*}
	J_{K,n}(e^{\hbar}) & = 1 + \lb 111 - 111 n^2 \rb \hbar^2 + \lb -2128 + 2128 n^2 \rb \hbar^3 + \lb \frac{253477}{4} - 73197 n^2 + \frac{39311}{4} n^4 \rb \hbar^4\\
	& + \lb - 2159616 + \frac{7947776}{3} n^2 - \frac{1468928}{3} n^4 \rb \hbar^5 + \cdots
\end{align*}
\end{prop}

\begin{lemma}
The Alexander polynomial of $C_{(16,3)}(\bm{4_1})$ is given by
\begin{equation*}
\begin{aligned}		
\Delta_{C_{(16,3)}{(\bm{4_1})}}(x) & = t^{21}-\frac{1}{t^{21}}+t^{20}+\frac{1}{t^{20}}+2 t^{18}+\frac{2}{t^{18}}-2t^{17}-\frac{2}{t^{17}}+t^{15}+\frac{1}{t^{15}}-t^{14}-\frac{1}{t^{14}}+t^{12}\\
& +\frac{1}{t^{12}} -t^{11}-\frac{1}{t^{11}}+t^9+\frac{1}{t^9}-t^8-\frac{1}{t^8}+t^6+\frac{1}{t^6}-t^5-\frac{1}{t^5}+t^3+\frac{1}{t^3}-2 t^2-\frac{2}{t^2}\\
& +t+\frac{1}{t}+1\\
\end{aligned}
\end{equation*}
\end{lemma}

\begin{prop}
The symmetric expansion of the inverse of the Alexander polynomial of $K$ in the limit of $\hbar \rarw 0$ is given by
\begin{align*}
\lim_{q \rarw 1} 2F_{K}(x,q) &  =  2\, \text{s.e.} \lb \frac{x^{1/2}-x^{-1/2}}{\Delta_{K}(x)} \rb\\
& = 2 x^{35/2} -\frac{2}{x^{35/2}} +4 x^{41/2} -\frac{4}{x^{41/2}} +10 x^{47/2} -\frac{10}{x^{47/2}} +26 x^{53/2} -\frac{26}{x^{53/2}} +68 x^{59/2}\\
& -\frac{68}{x^{59/2}}  +178 x^{65/2} -\frac{178}{x^{65/2}} -2 x^{67/2} +\frac{2}{x^{67/2}} +466 x^{71/2} -\frac{466}{x^{71/2}} -4 x^{73/2} +\frac{4}{x^{73/2}} \\
& +1220x^{77/2} -\frac{1220}{x^{77/2}}-10 x^{79/2}+\frac{10}{x^{79/2}} +3194 x^{83/2}-\frac{3194}{x^{83/2}} + \cdots \in \intg \lsb \lsb x^{\pm 1/2} \rsb \rsb 
\end{align*}
\end{prop}
\noindent In this case, the coefficients of $x^{67/2}, x^{73/2}, x^{79/2}$ etc are negative, which coincide with the negative coefficient elements in the initial data in Conjecture 4.1.
\newline

\begin{prop}
The initial data for the recursion relation of $K$ consists of $I=\lac f_{j} (q) \in \intg[q^{\pm 1}] | j=1,\cdots ,387 \rac$. The nonzero elements in $I$ are obtained by setting $w=5$ and thus $n=17$ in Conjecture 4.1.
\end{prop}
\noindent We verified that $f_{389}$ and $f_{391}$ obtained from the recursion match with those from the applications of Conjecture 4.4. Furthermore, vanishing of $f_{393}$ is consistent with Conjecture 4.4. 
\newline

\noindent \textbf{Acknowledgments.} I would like to thank the referee for helpful suggestions.


%
%
%

%
\appendix
\section*{Appendix}
\addcontentsline{toc}{section}{Appendix}

\section{The colored Jones polynomial of a cable knot}

For (r,2)-cabling of a knot $K$ in $S^3$, its unnormalized $\mathfrak{sl}_2(\mathbb{C})$ colored Jones polynomial of is~\cite{Ru}
$$
\tilde{J}_{C_{(r,2)}(K),n}(t)= t^{-rs\left(n^2-1\right)} \sum_{w=-\frac{n-1}{2}}^{\frac{n-1}{2}} t^{4rw(sw+1)} \tilde{J}_{K,\,(2sw+1)}(t),\quad t^4 =q
$$
The 0-framed unknot $U$ normalization is
$$
\tilde{J}_{U,n}(t) = \frac{t^{2n} - t^{-2n}}{t^2 - t^{-2}}.
$$

\section{The recursion relation for the figure eight knot}

We record the recursion relation and its initial data for the figure eight knot found in \cite{GM}. We note that a closed form formula was found in \cite{P}.

\begin{equation*}
\begin{aligned}
h_1 (q) & = 1\\
h_3 (q) & = 2\\
h_5 (q) & = \frac{1}{q}+3+q \\
h_7 (q) & = \frac{2}{q^2}+\frac{2}{q}+5+2 q+2 q^2 \\
h_9 (q) & = \frac{1}{q^4}+\frac{3}{q^3}+\frac{4}{q^2}+\frac{5}{q}+8+5 q+4 q^2+3 q^3+q^4 \\
h_{11} (q) & = \frac{2}{q^6}+\frac{2}{q^5}+\frac{6}{q^4}+\frac{7}{q^3}+\frac{10}{q^2}+\frac{10}{q}+15+10 q+10 q^2+7 q^3+6 q^4+2 q^5+2 q^6 \\
h_{13} (q) & = \frac{1}{q^9}+\frac{3}{q^8}+\frac{4}{q^7}+\frac{7}{q^6}+\frac{11}{q^5}+\frac{15}{q^4}+\frac{18}{q^3}+\frac{21}{q^2}+\frac{23}{q}+27+23 q+21 q^2+18 q^3+15 q^4\\
& +11 q^5+7 q^6+4 q^7+3 q^8+q^9\\
\end{aligned}
\end{equation*}

\begin{align*}
h_{m+14} (q) & =  -\frac{q^{-\frac{m}{2}-\frac{11}{2}}}{q^{\frac{m}{2}+\frac{13}{2}}-1} \lsb  h_m \left(q^{\frac{m}{2}+\frac{17}{2}}-q^{m+9}\right)+ h_{m+2} \left(q^{\frac{m}{2}+\frac{15}{2}}-q^{\frac{m}{2}+\frac{17}{2}}+q^{m+9}-q^{m+10}\right)  \notag\right.\\
& + h_{m+4} \left(-q^{\frac{m}{2}+\frac{11}{2}}-q^{\frac{m}{2}+\frac{17}{2}}-q^{\frac{m}{2}+\frac{19}{2}}+q^{\frac{3m}{2}+\frac{21}{2}}+q^{m+8}+q^{m+9}+q^{m+12}-q^7\right)\\
& + h_{m+6} \left(-q^{\frac{m}{2}+\frac{9}{2}}+q^{\frac{m}{2}+\frac{11}{2}}-q^{\frac{m}{2}+\frac{15}{2}}-q^{\frac{m}{2}+\frac{17}{2}} +q^{\frac{3 m}{2}+\frac{25}{2}}+q^{m+9}+q^{m+10}-q^{m+12}+q^{m+13}-q^5\right)\\
& + h_{m+8} \left(q^{\frac{m}{2}+\frac{11}{2}}+q^{\frac{m}{2}+\frac{13}{2}}-q^{\frac{m}{2}+\frac{17}{2}}+q^{\frac{m}{2}+\frac{19}{2}} -q^{\frac{3 m}{2}+\frac{31}{2}}-q^{m+8}+q^{m+9}-q^{m+11}-q^{m+12}+q^2\right)\\
& + h_{m+10} \left(q^{\frac{m}{2}+\frac{9}{2}}+q^{\frac{m}{2}+\frac{11}{2}}+q^{\frac{m}{2}+\frac{17}{2}}-q^{\frac{3m}{2}+\frac{35}{2}}-q^{m+9}-q^{m+12}-q^{m+13}+1\right)\\
&  + h_{m+12} \left(q^{\frac{m}{2}+\frac{11}{2}}-q^{\frac{m}{2}+\frac{13}{2}}+q^{m+11}-q^{m+12}\right)\notag\left.  \rsb
\end{align*}

\section{The definitions of the operators}

We list the definitions of the operators in the $\hat{A}$-polynomial (4) in Section 2.

$$
Q(t,M,L)= Q_2(t,M)L^2 + Q_1(t,M)L + Q_0(t,M),\quad B(t,M):= \sum_{j=0}^{2}c_{j}b(t,t^{2j+2}M^2)
$$
$$
b(t,M)=\frac{M(1+t^4 M^2)(-1+t^4 M^4)(-t^2 + t^{14} M^4)}{t^2 - t^{-2}}
$$
$$
c_{0}=\hat{P}_0 (t, t^4 M^2)\hat{P}_1 (t,t^6 M^2),\quad c_{1}= - \hat{P}_1(t,t^2 M^2)\hat{P}_{1}(t,t^6 M^2),\quad c_{2}=\hat{P}_1(t,t^2 M^2)\hat{P}_2(t,t^4 M^2).
$$
\begin{equation*}
\begin{aligned}
Q_2(t,M) & = \, \hat{P}_2(t,t^4 M^2)\, \hat{P}_1(t,t^2 M^2)\, \hat{P}_0(t,t^6 M^2)\\
Q_1(t,M) & = \, \hat{P}_0(t,t^4 M^2)\, \hat{P}_1(t,t^6 M^2)\, \hat{P}_2(t,t^2 M^2)- \hat{P}_1(t, t^6 M^2)\, \hat{P}_1(t,t^2 M^2)\, \hat{P}_1(t,t^4 M^2)\\
& + \hat{P}_2(t,t^4 M^2)\, \hat{P}_1(t,t^2 M^2)\, \hat{P}_0(t,t^6 M^2)\\
Q_0(t,M) & = \, \hat{P}_0(t,t^4 M^2)\, \hat{P}_1(t, t^6 M^2)\, \hat{P}_0(t,t^2 M^2),
\end{aligned}
\end{equation*}
\begin{equation*}
\begin{aligned}
\hat{P}_0(t,M): & = \, t^6 M^4 (-1 + t^{12} M^4)\\
\hat{P}_1(t,M): & = \, -(-1+ t^4 M^2)(1+ t^4 M^2) \lb 1-t^4 M^2 - t^4 M^4 - t^{12}M^4- t^{12}M^4 - t^{12}M^6 + t^{16}M^8 \rb\\
\hat{P}_2(t,M): & = \, t^{10}M^4 ( -1 + t^4 M^4 ).
\end{aligned}
\end{equation*}
\newline

The definitions of the operators in the $\hat{A}$-polynomial (12) in Section 4 are

$$
Q(t,M,L) : = 1 + \frac{Det A_9}{Det A} L + \frac{Det A_8}{Det A} L^2, \quad B(t,M):= \sum_{j=0}^{6}c_{j}b(t,t^{2(2+j)}M^3)
$$
where $A$ is a matrix and $A_9$ and $A_8$ are obtained from $A$ (see Section 4 in \cite{Ru}).

%
%

\end{document}